\def\A{\leavevmode\setbox0\hbox{A}\lower1.4ex\hbox 
to\wd0{\hss`}\kern-.9\wd0A}
\def\E{\leavevmode\setbox0\hbox{E}\lower1.4ex\hbox 
to\wd0{\hss`\/}\kern-.9\wd0E}
\def\a{\leavevmode\setbox0\hbox{a}\lower1.4ex\hbox 
to\wd0{\hss`\/}\kern-\wd0a}
\def\e{\leavevmode\setbox0\hbox{e}\lower1.4ex\hbox 
to\wd0{\hss`\/}\kern-\wd0e}
\newcommand{\be}{\begin{equation}}
\newcommand{\ee}{\end{equation}}
\newcommand{\ba}{\begin{array}}
\newcommand{\ea}{\end{array}}
\newcommand{\beqn}{\begin{eqnarray}}
\newcommand{\eeqn}{\end{eqnarray}}
\newcommand{\nn}{\nonumber}
\newcommand{\p}{\prime}
\newcommand{\x}{(m^{\p\p},\tau^{\p\p},a^{\p\p},v^{\p\p})}
\newcommand{\y}{(m^{\p},\tau^{\p},a^{\p},v^{\p})}
\newcommand{\z}{(m,\tau,a,v)}
\newcommand{\pa}{\partial}

\font\symb=msam7
\def\znakr{\raise1.5pt\hbox{\symb\char66\kern-2pt\char74}}
\def\znakl{\raise1.5pt\hbox{\symb\char73\kern-2pt\char67}}

\def\normalsize{
\setlength{\textheight}{23cm}
\setlength{\textwidth}{15cm}
\setlength{\topmargin}{-2.0cm}
\setlength{\hoffset}{-0.5cm}
\setlength{\leftmargin}{-1cm}
\setlength{\rightmargin}{2.0cm}}

\documentstyle[12pt]{article}
\normalsize

\begin{document}
\title{Lie Bialgebra Structures for Centrally Extended
Two-Dimensional Galilei Algebra and their Lie--Poisson
Counterparts}
\author{ Anna Opanowicz\thanks{Supported by
\L\'od\'z University Grant No.\ 580} \\
Department of Theoretical Physics II \\
University of \L \'od\'z \\
Pomorska 149/153, 90-236 \L \'od\'z, Poland}
\date{}
\maketitle
\setcounter{section}{0}
\setcounter{page}{1}
\begin{abstract}
All bialgebra structures for centrally extended Galilei algebra
are classified. The corresponding Lie--Poisson structures on
centrally extended Galilei group are found.
\end{abstract}

\newpage
\section{Introduction}

Much interest has been attracted in last years to the problem
of deformations of space-time symmetry groups \cite{1},\cite{2},\cite{3},\cite{4},\cite{5},\cite{6},\cite{7}.
In particular, in the recent paper \cite{8} all inequivalent bialgebra
structures on two-dimensional Galilei algebra were classified and the
corresponding Lie--Poisson structures on the group were found.

From the physical point of view what is really interesting is
the central extension of the Galilei algebra. This is because
only the genuine projective representations of the Galilei
group are relevant in nonrelativistic quantum theory \cite{9}.
In the present paper we classify all nonequivalent bialgebra
Lie--Poisson structures for centrally extended two-dimensional
Galilei algebra/group. In the two-dimensional case there exists
two-parameter family of central extensions, the parameters
being the mass of the particle and the constant force acting
on it. We restrict ourselves to the case of free particles, i.e.
only the mass parameter is kept nonvanishing.

The content of the paper is as follows. First, we find the general
form of 1-cocycle on centrally extended twodimesional Galilei
algebra. Then the action of most general automorphism
transformation on such 1-cocycle is considered and its orbits
are classified which allows to find all nonequivalent bialgebra
structures. The corresponding Lie--Poisson structures on Galilei
group are then found. The whole procedure follows quite closely
the one presented in Ref.~[10] for $E(2)$ group and in Ref.~[8]
for  two-dimensional Galilei group. As a result we find 26
nonequivalent bialgebra structures (some of them being still
one parameter families), 8 of them being the coboundary ones.

\section{Two-dimensional Galilei group and algebra with central extension
and their automorphisms}

The two-dimensional Galilei group is a Lie group of transformations
of the space-time with one space dimensions. An arbitrary group
element $g$ is of the form
\be
g=(\tau , v, a);
\label{r1}
\ee
here $\tau $ is time translation, $a$ respectively $v$ are space translation
respectively Galilean boost.
The multiplication law reads:
\be
g^\prime  g=(\tau' +\tau , v^\prime  +v, a^\prime  +a+\tau
v^\prime  );
\label{r2}
\ee
The resulting Lie algebra takes the form:
\be
[K,H]=iP,\ \  [K,P]=0,\ \ [H,P]=0
\label{r3}
\ee

The central extension is obtained by replacing the second commutation
rule by
\be
[K,P]=iM
\label{r4}
\ee
where
\be
[M,\cdot\;]=0.
\label{r5}
\ee
Therefore, we arrive finally at the following algebra
\be
[K,H]=iP,\ \ [K,P]=iM,\ \ [H,P]=0,\ \ [M,\cdot\;]=0
\label{r6}
\ee
Let us define the centrally extended Galilei group by the following
global expotential parametrization of group elements
\be
\widetilde g\ = e^{imM} e^{-i\tau H} e^{iaP} e^{ivK}
\label{r7}
\ee
Let us write
\be
\widetilde g\ = (m, \tau , v, a)
\label{r8}
\ee
Then we have the following multiplication law
\be
\widetilde{g}^\prime
\widetilde{g} =
(m^\prime + m -{1\over 2} v^{\prime 2} \tau -a v^\prime ,
\tau ^\prime  + \tau , v^\prime  + v, a^\prime  + a + \tau v^\prime )
\label{r9}
\ee

Lie algebra with central extension can be realized in terms of right -
invariant fields to be calculated according to the standard rules from
the composition law (\ref{r9})
\beqn
&&X^{R}_{v} = i ( {\partial\over{\partial v}} -
a {\partial\over{\partial m}}
+\tau  {\partial\over{\partial a}} ) \nonumber\\
&&X^{R}_{a} = i {\partial\over{\partial a}} \nonumber\\
							\label{r10}
&&X^{R}_{m} = i {\partial\over{\partial m}} \nonumber\\
&&X{^{R}_{\tau }} = - i {\partial\over{\partial \tau }}
\eeqn

Let us now describe all automorphisms of the algebra (\ref{r6}).
The group of automorphisms consists of the following transformations
\be
\left(
\ba{lcr}
K\\
H\\
P\\
M
\ea
\right) \rightarrow
\left(
\ba{lcr}
K^\prime \\
H^\prime \\
P^\prime \\
M^\prime
\ea
\right)
=
\left(
\ba{lccccr}
\gamma_3 & \alpha_3 & \beta_3 & \eta_3 \\
0 & \alpha_1 & \beta_1 & \eta_1 \\
0 & 0 & \beta_2 & \eta_2 \\
0 & 0 & 0 & \eta_4
\ea
\right)
\left(
\ba{lcr}
K\\
H\\
P\\
M
\ea
\right)
\label{r11}
\ee
where :
\beqn
&&\beta_2 = \gamma_3 \alpha_1\nonumber\\
&&\eta_2 = \gamma_{3} \beta_{1}\label{r12}\\
&&\eta_4 = \gamma_3 \beta_2,\nonumber
\eeqn
and, obviously, \( \alpha_1 \neq 0, \gamma_3 \neq 0 \).

\section{\sloppy The Bialgebra structures on two-di\-men\-sional
centrally extended Galilei algebra}

Our aim here is to give a complete classification of Lie bialgebra structures
for the algebra (\ref{r6}) up to automorphisms.

Let us remind the definition of bialgebra. It is a pair \( (L, \delta)\),
where \(L\) is a Lie algebra while  \(\delta\)   is a skewsymmetric
cocommutator
\( \delta : L \rightarrow L \otimes L \),  i.e.\\
(i)  \( \delta \)   is a 1-cocycle,
\[ \delta([X,Y]) = [\delta(X), 1 \otimes Y + Y \otimes 1] +
[1 \otimes X + X \otimes 1, \delta(Y)]\; \; \rm{for} \; \; X,Y\in L  \]
(ii) The dual map \( \delta ^\ast : L^\ast \otimes L^\ast \rightarrow L^\ast\
\)
defines a Lie bracket on \( L^\ast \).

We can find all bialgebra structures on our algebra. The general form of
\(\delta\)
obeying (i) is
\be
\left(
\ba{lcr}
\delta(H)\\
\delta(P)\\
\delta(K)\\
\delta(M)
\ea
\right)
=
\left(
\ba{lccccccr}
a&0&0&b&c&d\\
0&0&0&a&0&h-b\\
e&0&f&g&h&j\\
0&0&0&0&0&-(a+f)
\ea
\right)
\left(
\ba{lcccr}
H \wedge P\\
H \wedge K\\
P \wedge K\\
H \wedge M\\
M \wedge K\\
M \wedge P
\ea
\right),
\label{r13}
\ee
a, b, c, d, e, f, g, h and j being arbitrary real parameters.

From the condition (ii) we obtain:
\[ a = b = c = 0\]
or
\be a = e = f = 0 \label{r14} \ee
or
\[ b = c = f = h = 0\]

Eqs.~(\ref{r13}) and (\ref{r14}) define all bialgebra structures on
two-dimensional Galilei algebra (\ref{r6}). However, we are interested
in classification of nonequivalent bialgebra structures. To this end
we find the transformation rules for the parameters under the
automorphisms (\ref{r11}). They read
\be
\ba{l}
\widetilde a ={{a} \over {\beta_2}}\nonumber\\
\widetilde b = {{b\over \eta_4} +{{\alpha_3\alpha_1}\over
{\beta_2\eta_4}}c}\nn\\
\widetilde c = {\alpha_1\over {\eta_4\gamma_3}}c\nn\\
\widetilde d = {{{\alpha_1\over{\beta_2\eta_4}}d}-2{{{\eta_1\alpha_1}\over
{{\beta_2}^3}}a}+{{\alpha_1(\alpha_3\beta_1-\alpha_1\beta_3)}
\over{{{\beta_2}^2}\eta_4}}c}+{{{\beta_1\eta_2}\over{{{\beta_2}^2}\eta_4}}a}
+ {{\beta_1\over{\beta_2\eta_4}}h}-{{\eta_1\over{\beta_2\eta_4}}f}\nn\\
\widetilde e = {{{\gamma_3\over{\alpha_1\beta_2}}e}+{{{\alpha_3\gamma_3}
\over{{\beta_2}^2}}f}+{{\alpha_3\over{\alpha_1\beta_2}}a}}\nn\\
\widetilde f = {f\over{\beta_2}}\label{r15}\\
\widetilde g = {{{{\gamma_3\over{\alpha_1\eta_4}}g} - {{\beta_1\gamma_3}\over
{\alpha_1{{\beta_2}^2}}}e} - {{\alpha_3\eta_2}\over{{\beta_2}^3}}f} +
{{\alpha_3\over{\alpha_1\eta_4}}b} + {{{\alpha_3\gamma_3}
\over{\beta_2\eta_4}}h}-{{{\beta_1\alpha_3}\over{{{\beta_2}^2}\alpha_1}}a}
+ {{{{\alpha_3}^2}\over{\beta_2\eta_4}}c}+{{{\beta_3}\over{\alpha_1\eta_4}}a}
\nn \\
\widetilde j = {{\gamma_3\over{\beta_2\eta_4}}j - {{\eta_1\gamma_3}\over
{{\beta_2}^3}}e-{{\alpha_3\eta_1}\over{{\beta_2}^3}}f +
{{\gamma_3\eta_2}\over
{{{\beta_2}^2}\eta_4}}g +{{\alpha_3\beta_1}\over{{\beta_2}^3}}h
-{{\alpha_3\eta_1}\over{{\beta_2}^3}}a+{{\alpha_3\eta_2}
\over{{{\beta_2}^2}\eta_4}}b}\nonumber\\
\; \; + {{{\alpha_3(\alpha_3\beta_1-\alpha_1\beta_3)}
\over{{{\beta_2}^2}\eta_4}}c+
{{\alpha_3\over{\beta_2\eta_4}}d+{{\eta_2\beta_3}\over{{{\beta_2}^2}\eta_4}}a-
{\beta_3\over{\beta_2\eta_4}}b-{\eta_3\over{\beta_2\eta_4}}}a}\nn\\
\widetilde h = {{h\over\eta_4}-{\beta_1\over{{\beta_2}^2}}f
+{\alpha_3\over{\eta_4\gamma_3}}c}\nn\\
\widetilde h-\widetilde b = {{{h-b}\over\eta_4}-
{\eta_2\over{\beta_2\eta_4}}f}\nn\\
\widetilde a+\widetilde f = {{a+f}\over\beta_2}\nn
\ea
\ee
We are now in position to classify all orbits of automorphism group
in the space of bialgebra structures.
A simple but long and painful analysis leads to the complete list of Lie
bialgebra structures summarized in Table 1.\\
\newpage
\centerline{Table 1}
\begin{center}
\begin{tabular}{|c|c|c|c|c|c|c|c|c|c|c|}
\hline
 & a & b & c & d & e & f & g & h & j & Remarks \\  \hline \hline
1 & 0 & 0 & 0 & 0 & 1 & 0 & 0 & 0 & 0 &  \\ \hline
2 & 0 & 0 & 0 & 0 & -1 & 0 & 0 & 0 & 0 &  \\ \hline
3 & 0 & 0 & 0 & 0 & 0 & 0 & 0 & 1 & 0 & coboundary \\ \hline
4 & 0 & 0 & 0 & 0 & 0 & 0 & 1 & 0 & 0 & coboundary \\  \hline
5&0&0&0&1&0&0&0&0&0&coboundary\\ \hline
6&0&0&0&0&0&0&0&0&1&coboundary\\ \hline
7&0&0&0&0&0&0&0&0&-1&coboundary\\ \hline
8&0&0&0&0&0&0&0&1&1&coboundary\\ \hline
9&0&0&0&0&0&0&0&1&-1&coboundary\\ \hline
10&0&0&0&0&1&0&0&1&0&  \\ \hline
11&0&0&0&0&-1&0&0&1&0&  \\ \hline
12&0&0&0&1&0&0&1&0&0&coboundary\\ \hline
13&0&0&0&0&0&1&0&0&$\varepsilon$&$\varepsilon\in $ R\\ \hline
14&0&0&0&0&0&1&1&0&$\varepsilon$&$\varepsilon\in $ R\\ \hline
15&0&0&0&1&1&0&0&0&0&  \\ \hline
16&0&0&0&1&-1&0&0&0&0&  \\ \hline
17&0&1&0&0&0&0&0&0&0&  \\ \hline
18&0&1&0&1&0&0&0&0&0&  \\ \hline
19&0&$\varepsilon$&0&0&0&0&0&1&0&$\varepsilon \neq  0 $\\ \hline
20&0&-1&0&0&0&0&1&1&0&  \\ \hline
21&0&$\varepsilon$&1&0&0&0&0&0&0&$\varepsilon\in $ R\\ \hline
22&0&$\varepsilon$&1&0&0&0&0&0&1&$\varepsilon\in $ R\\ \hline
23&0&$\varepsilon$&1&0&0&0&0&0&-1&$\varepsilon\in $ R\\ \hline
24&0&$\varepsilon$&1&0&0&0&1&0&0&$\varepsilon\in $ R\\ \hline
25&0&$\varepsilon$&1&0&0&0&-1&0&0&$\varepsilon\in $ R \\ \hline
26&1&0&0&0&0&0&0&0&0&  \\ \hline
\end{tabular}
\end{center}

We have checked explictly that all the above bialgebra structures
are consistent and inequivalent. It remains to find coboundary
structures (listed also in Table1).

As it is well known a cocommutator $\delta$ given by
\be
\delta(X) = i[1\otimes X + X\otimes 1,\; r],\; r\in L\wedge L,\;  X\in L
\label{r16}
\ee
defines a coboundary Lie bialgebra if and only if $r$ fulfills the
modified classical Yang-Baxter equation
\be
[X\otimes 1\otimes 1 + 1\otimes X\otimes 1 +
1\otimes 1\otimes X, \; \xi(r)] = 0, \; X\in L,
\label{r17}
\ee
where $\xi(r)$ is the Schouten bracket
\[\xi (r) \equiv [r_{12}, \; r_{13}] + [r_{12}, \; r_{23}]
 + [r_{13}, \; r_{23}];\]
here
\beqn
r_{12} = r^{ij}{X_i}\otimes {X_j}\otimes 1\nn\\
r_{13} = {r^{ij}{X_i}\otimes 1\otimes {X_j}}\nn\\
r_{23} ={{r^{ij}}1\otimes {X_i}\otimes {X_j}}.\nn
\eeqn
Let us put
\be
r = A H\wedge P + C P\wedge K + D H\wedge M + E M\wedge K + F M\wedge P
\label{r18}
\ee
Egs. (\ref{r18}) and (\ref{r16}) give now
\beqn
&&\delta(P) = -C M\wedge P\nonumber\\
&&\delta(H) = E M\wedge P\nonumber\\
&&\delta(K) = -A H\wedge M - C M\wedge K + D M\wedge P
\label{r19} \\
&&\delta(M) = 0\nonumber
\eeqn
By comparying Eqs.~(\ref{r13}) and (\ref{r19}) we get
\beqn
&&a=b=c=e=f=0\nn\\
&&A=-g, \; C=-h, \; D=j, \; E=d\nn
\eeqn
which serve to identify the coboundary structures in Table 1.

\section{The Lie--Poisson structures on two-di\-men\-sional Galilei group}
In this section we find all Lie--Poisson structures on centrally
extended two-di\-men\-sional Galilei group. Let $G$ be a Lie group , $L$
its Lie algebra and \\
$\{X^{R}_{i}\}$- the set of right invariant fields
on $G$. As it is well known
\be
\{\Psi, \; \Phi \}\equiv {\eta ^{ij}(g)}{{X^{R}_{i}}\Psi}{{X^{R}_{j}}\Phi}
\label{r20}
\ee
where
\be
\eta(g) = \eta^{ij}(g)X_{i}\otimes X_{j}, \;\eta:G\rightarrow \Lambda^{2}L
\label{r21}
\ee
provides $G$ with a Poisson-Lie group structure if only if
\beqn
&&(i) \; \eta^{il}X^{R}_{l}\eta^{jk} + \eta^{kl}X^{R}_{l}\eta^{ij} +
\eta^{jl}X^{R}_{l}\eta^{ki} \nn\\
&&\; \; \; \; - c^{j}_{lp}\eta^{il}\eta^{pk} - c^{i}_{lp}\eta^{kl}\eta^{pj}
- c^{k}_{lp}\eta^{jl}\eta^{pi} = 0\label{r22} \\
&&(ii) \; \eta(g^{\prime}g) = \eta(g^{\prime}) + Adg^{\prime}\eta(g)
\label{r23}
\eeqn
In our case let us define
\beqn
\eta(m,\tau,a,v) &=& \lambda(m,\tau,a,v)H\wedge P +
\mu(m,\tau,a,v)H\wedge K \nn \\
&& + \nu(m,\tau,a,v)P\wedge K + \kappa(m,\tau,a,v)H\wedge M
\label{r24}\\
&& + \rho(m,\tau,a,v)M\wedge K
+ \pi(m,\tau,a,v)M\wedge P \nonumber
\eeqn
Eq.~(\ref{r23}) gives
\begin{eqnarray}
\eta(g^{\prime}g) &=& (\lambda^{\prime} + \lambda
-\tau^{\prime}\mu)H\wedge P +
(\nu^{\prime} + \nu - v^{\prime}\mu)P\wedge K \nonumber \\
&&+   (\kappa^{\prime} + \kappa - v^{\prime}\lambda
+ a^{\prime}\mu)H\wedge M + (\mu^{\prime} + \mu)H\wedge K\nonumber\\
&&+ (\rho^{\prime} + \rho +
\frac{{v'}^2}{2}\mu- v^{\prime}\nu)M\wedge K \label{r25} \\
&& + (\pi^{\prime} + \pi -
\frac{{v^\p}^2}{2}\lambda
+ (v^{\p}a^{\p} - \frac{{v^{\p}}^2\tau^{\p}}{2})\mu + (v^{\p}\tau^{\p}
- a^\p)\nu + v^{\p}\kappa - \tau^{\p}\rho)M\wedge P \nn
\end{eqnarray}
where $\lambda^\p \equiv \lambda\y$ etc.

Consequently we obtain the following set of equations determining $\lambda$,
$\mu$, $\nu$, $\kappa$, $\rho$ and $\pi$
\beqn
\lambda\x &=& \lambda\y + \lambda\z - {\tau^\prime}\mu\z \nn \\
\mu\x &=& \mu\y + \mu\z \nn \\
\nu\x &=& \nu\y + \nu\z - {v^\p}\mu\z \label{r26} \\
\kappa\x &=& \kappa\y + \kappa\z -{v^\p}\lambda\z \nn\\
&&+ {a^\p}\mu\z \nn\\
\rho\x &=& \rho\y + \rho\z + {{1\over2}{v^\p}^2}\mu\z\nn\\
&&- v'\nu\z\nn\\
\pi\x &=& \pi\y + \pi\z - {{1\over 2}{{v^\p}^2}}\lambda\z\nn\\
&&+ {v^\p}\kappa\z + ({v^\p}{a^\p}-{1\over2}{{v^\p}^2}{\tau^\p})\mu\z\nn\\
&& + ({v^\p}{\tau^\p}
- {a^\p})\nu\z - {\tau^\p}\rho\z \nn
\eeqn

The strategy to solve Eq.~(\ref{r26}) is to find first the form of $\eta$
for 1-parameter subgroups generated by $P$, $K$, $H$, $M$ and use again
Eq.~(\ref{r26}) together with the decomposition
\be
\z = (m,0,0,0)\ast (0,\tau,0,0)\ast (0,0,a,0)\ast (0,0,0,v)
\label{r27}
\ee
to determine the form of $\eta$ for general group element. Eq.~(\ref{r26})
as specialized for one-parameter subgroups generated by $M$, $H$, $P$ and $K$
read
\beqn
&&\lambda(m^{\p\p},0,0,0) = \lambda(m^{\p},0,0,0) + \lambda(m,0,0,0)\nn \\
&&\mu(m^{\p\p},0,0,0)  = \mu(m^{\p},0,0,0)  + \mu(m,0,0,0)\nn \\
&&\nu(m^{\p\p},0,0,0)  = \nu(m^{\p},0,0,0)  + \nu(m,0,0,0)\nn \\
&&\kappa(m^{\p\p},0,0,0)  = \kappa(m^{\p},0,0,0)  + \kappa(m,0,0,0)\nn\\
&&\rho(m^{\p\p},0,0,0)  = \rho(m^{\p},0,0,0)  + \rho(m,0,0,0)\nn \\
&&\pi(m^{\p\p},0,0,0)  = \pi(m^{\p},0,0,0)  + \pi(m,0,0,0)\nn \\
&&\nn\\
&&\lambda(0,\tau^{\p\p},0,0)  = \lambda(0,\tau^{\p},0,0) + \lambda(0,\tau,0,0)
- {\tau^\p}\mu(0,\tau,0,0) \nn \\
&&\mu(0,\tau ^{\p\p},0,0) = \mu(0,\tau ^{\p},0,0)  + \mu(0, \tau ,0,0)\nn \\
&&\nu(0,\tau ^{\p\p},0,0)  = \nu(0,\tau ^{\p},0,0)  + \nu(0, \tau ,0,0)\nn \\
&&\kappa(0,\tau ^{\p\p},0,0)  = \kappa(0,\tau ^{\p},0,0)  +
\kappa(\tau ,0,0,0)\nn\\
&&\rho(0,\tau ^{\p\p},0,0) = \rho(0,\tau ^{\p},0,0)  +
\rho(0, \tau ,0,0)\nn \\
&&\pi(0,\tau ^{\p\p},0,0)  = \pi(0,\tau ^{\p},0,0)  + \pi(0, \tau ,0,0)
- {\tau^\p}\rho(0,\tau,0,0)  \label{r28} \\
&&\nn\\
&&\lambda(0,0,a^{\p\p},0)  = \lambda(0,0,a^{\p},0) + \lambda(0,0,a,0)\nn \\
&&\mu(0,0,a^{\p\p},0)  = \mu(0,0,a^{\p},0)  + \mu(0,0,a,0)\nn \\
&&\nu(0,0,a^{\p\p},0)  = \nu(0,0,a^{\p},0)  + \nu(0,0,a,0)\nn \\
&&\kappa(0,0,a^{\p\p},0)  = \kappa(0,0,a^{\p},0)  + \kappa(0,0,a,0)
+ {a^\p}\mu(0,0,a,0)\nn\\
&&\rho(0,0,a^{\p\p},0)  = \rho(0,0,a^{\p},0)  + \rho(0,0,a,0)\nn \\
&&\pi(0,0,a^{\p\p},0)  = \pi(0,0,a^{\p},0)  + \pi(0,0,a,0)
- {a^\p}\nu(0,0,a,0)\nn \\
&&\nn\\
&&\lambda(0,0,0,v^{\p\p})  = \lambda(0,0,0,v^{\p}) + \lambda (0,0,0,v)\nn \\
&&\mu(0,0,0,v^{\p\p})  = \mu(0,0,0,v^{\p})  + \mu(0,0,0,v)\nn \\
&&\nu(0,0,0,v^{\p\p})  = \nu(0,0,0,v^{\p})  + \nu(0,0,0,v)
- {v^\p}\mu(0,0,0,v) \nn \\
&&\kappa(0,0,0,v^{\p\p})  = \kappa(0,0,0,v^{\p})  + \kappa(0,0,0,v)
- {v^\p}\lambda(0,0,0,v) \nn\\
&&\rho(0,0,0,v^{\p\p})  = \rho(0,0,0,v^{\p})  + \rho(0,0,0,v)
+ {1\over 2}{{v^\p}^2}\mu(0,0,0,v) - {v^\p}\nu(0,0,0,v)\nn \\
&&\pi(0,0,0,v^{\p\p})  = \pi(0,0,0,v^{\p})  + \pi(0,0,0,v)
- {1\over 2}{{v^\p}^2}\lambda(0,0,0,v) + {v^\p}\kappa(0,0,0,v)\nn
\eeqn
The corresponding solutions read
\beqn
&&\lambda(m,0,0,0) = {a_1}m\nn\\
&&\mu(m,0,0,0) = {a_2}m\nn\\
&&\nu(m,0,0,0) = {a_3}m\nn\\
&&\kappa(m,0,0,0) ={a_4}m\nn\\
&&\rho(m,0,0,0) = {a_5}m\nn\\
&&\pi(m,0,0,0) = {a_6}m\nn\\
&&\nn\\
&&\lambda(0,\tau,0,0) = {b_1}\tau - {1\over 2}{b_2}{{\tau }^2}\nn\\
&&\mu(0,\tau,0,0) = {b_2}\tau\nn\\
&&\nu(0,\tau,0,0) = {b_3}\tau \nn\\
&&\kappa(0,\tau,0,0) = {b_4}\tau\nn\\
&&\rho(0,\tau,0,0) ={b_5}\tau\nn\\
&&\pi(0,\tau,0,0) = {b_6}\tau - {1\over 2}{b_5}{{\tau }^2}\label{r29}\\
&&\nn\\
&&\lambda(0,0,a,0) = {c_1}a\nn\\
&&\mu(0,0,a,0) = {c_2}a\nn\\
&&\nu(0,0,a,0) = {c_3}a \nn\\
&&\kappa(0,0,a,0) = {c_4}a+ {1\over 2}{c_2}{a^2}\nn\\
&&\rho(0,0,a,0) = {c_5}a\nn\\
&&\pi(0,0,a,0) = {c_6}a - {1\over 2}{c_3}{{a}^2}\nn\\
&&\nn\\
&&\lambda(0,0,0,v) = {d_1}v\nn\\
&&\mu(0,0,0,v) = {d_2}v\nn\\
&&\nu(0,0,0,v) = {d_3}v - {1\over 2}{d_2}{v^2}\nn\\
&&\kappa(0,0,0,v) = {d_4}v- {1\over 2}{d_1}{v^2}\nn\\
&&\rho(0,0,0,v) = {d_5}v + {1\over 6}{d_2}{v^3} - {1\over 2}{d_3}{v^2}\nn\\
&&\pi(0,0,0,v) = {d_6}v - {1\over 6}{d_1}{v^3} + {1\over 2}{d_4}{v^2}.\nn
\eeqn

Now, using Eqs.~(\ref{r29}) and (\ref{r26}) we obtain the general
form of $\lambda,\mu,\nu,\kappa,\rho$ and $\pi$ (the resulting
expressions were reinserted back to Eq.~(\ref{r26}) which provided
some furhter constraints for parameters):
\beqn
&&\lambda\z = {b_1}\tau + {d_1}v\nn\\
&&\mu\z = 0\nn\\
&&\nu\z = {d_3}v\label{r30}\\
&&\kappa\z = {b_4}\tau - {b_1}a + {d_4}v - {1\over 2}{d_1}{v^2}  \nn \\
&&\rho\z = {b_5}\tau + {d_5}v - {1\over 2}{d_3}{v^2}\nn\\
&&\pi\z = ({b_1}-{d_3})m + {b_6}\tau - {1\over 2}{b_5}{\tau^2} + ({b_4}
+ {d_5})a + {d_6}v - {1\over 6}{d_1}{v^3}\nn\\
&&\;\;\;\;\; + {1\over 2}{d_4}{v^2} - {d_3}av - {d_5}v\tau + {1\over 2}
{d_3}\tau{v^2}\nn
\eeqn

The general form of $\eta$ is given by Eqs.~(\ref{r24}) and (\ref{r30}).
Our next aim is to classify nonequivalent {$\eta^\p$}s. As it is well
known $\eta $ defines the bialgebra structure on $L$  through
\be
\delta(X)={{d\eta(e^{itx})}\over{dt}}\mid _{t=0}
\label{r31}
\ee
Simple calculation gives:
\newpage
\beqn
&&\delta(H) = -{b_1}H\wedge P - {b_4}H\wedge M - {b_5}M\wedge K
- {b_6}M\wedge P \nn \\
&&\delta(P) = -{b_1}H\wedge M + ({b_4} + {d_5})M\wedge P \label{r32}\\
&&\delta(K) = {d_1}H\wedge P + {d_3}P\wedge K + {d_4}H\wedge M
+ {d_5}M\wedge K + {d_6}M\wedge P \nn \\
&&\delta(M) = ({b_1} - {d_3})M\wedge P \nn
\eeqn
By comparying Eqs.~(\ref{r13}) and (\ref{r32}) we get
\beqn
&&a = -{b_1}, \; e = {d_1}, \; {j = {d_6}}\nn\\
&&b = -{b_4}, \; f = {d_3}, \; -(a+f) = {{b_1} - {d_3}}\label{r33}\\
&&c = -{b_5}, \; g = {d_4}, \; h-b = {{b_4} + {d_5}}\nn\\
&&d = -{b_6}, \; h = {d_5}, \nn
\eeqn

Eq.~(\ref{r33}), together with the results of previous section (Table 1)
gives us all inequivalent Poisson structures on two-dimensional centrally
extended Galilei group. To this end we write out explicitly the general
form of Poisson bracket following from Eqs.~(\ref{r10}), (\ref{r20})
and (\ref{r30})
\beqn
&&\{f,g\} = \lambda({{\pa f}\over{\pa \tau}}{{\pa g}\over {\pa a}}
- {{\pa f}\over {\pa a}}{{\pa g}\over {\pa \tau}})
+ \kappa ({{\pa f}\over {\pa\tau }}{{\pa g}
\over {\pa m}} - {{\pa f}\over {\pa m}}{{\pa g}\over {\pa \tau }})\nn\\
&&\; \; \; \; + \mu\left( {{\pa f}\over
{\pa \tau }}({{\pa g}\over {\pa v}} - a{{\pa g}\over {\pa m}}
+ \tau {{\pa g}\over {\pa a}}) - ({{\pa f}\over {\pa v}}
-a{{\pa f}\over {\pa m}} + \tau {{\pa f}\over {\pa a}}){{\pa g}\over
{\pa \tau }}\right) \nn\\
&&\; \; \; \;  - \nu\left( {{\pa f}\over {\pa a}}({{\pa g}\over {\pa v}}
- a{{\pa g}\over {\pa m}} + \tau {{\pa g}\over {\pa a}}) -({{\pa f}
\over {\pa v}} - a{{\pa f}\over {\pa m}} + \tau {{\pa f}\over {\pa a}})
{{\pa g}\over {\pa a}}\right)
\label{r34}\\
&&\; \; \; \; -\rho\left( {{\pa f}\over {\pa m}}({{\pa g}\over {\pa v}}
-a{{\pa g}
\over {\pa m}} + \tau {{\pa g}\over {\pa a}}) - ({{\pa f}\over {\pa v}}
- a{{\pa f}\over {\pa m}} + \tau {{\pa f}\over {\pa a}})
{{\pa g}\over {\pa m}}\right)\nn\\
&&\; \; \; \; - \pi({{\pa f}\over {\pa m}}{{\pa g}\over{\pa a}}
- {{\pa f}\over {\pa a}}{{\pa g}\over {\pa m}})\nn
\eeqn
In particular, the basic Lie--Poisson brackets read:
\beqn
&&\{v,\tau\} = -\mu = 0\nn\\
&&\{v,a\} = \nu = {d_3}v\nn\\
&&\{v,m\} = \rho = {b_5}\tau + {d_5}v - {{1\over 2}{d_3}{v^2}}\nn\\
&&\{\tau,a\} = \lambda + \tau \mu = {b_1}\tau + {d_1}v \label{r35}\\
&&\{\tau,m\} = -a\mu + \kappa = {b_4}\tau - {b_1}a + {d_4}v -
{{1\over 2}{d_1}{v^2}}\nn\\
&&\{a,m\} = \pi + a\nu + \tau \rho = {{1\over 2}{b_5}{\tau ^2}}
+ ({b_1} - {d_3})m + {b_6}\tau  + ({b_4} + {d_5})a \nn\\
&&\; \; \; \; \; + {d_6}v - {{1\over 6}{d_1}{v^3}} + {{1\over 2}{d_4}{v^2}}\nn
\eeqn

Obviously, there are further constraints on parameters following from
Eq.~(\ref{r22}) which haven$^{\p}$t been used yet. Instead of solving
it we impose the Jacobi identities on our Poisson brackets (which is
equivalent to solving Eq.~(\ref{r22})). It appears that the additional
constaints are, through Eq.~(\ref{r33}), equivalent to the ones given
by Eq.~(\ref{r14}) which provides a further test of the consistency
of our results.

Eqs. (\ref{r33}), (\ref{r35}) and the classification given in Tabele 1
lead us finally to the following classification of nonequivalent Lie--Poisson
structures (Table 2)

\begin{center}
Table 2\\[2,5mm]
{\scriptsize
\begin{tabular}{|c|c|c|c|c|c|c|c|}
\hline
 &  \{ $v$, $a$\}& \{ $v$, $m$\}& \{$\tau$, $a$\}&
\{ $\tau$, $m$\}&\{ $a$, $m$\} &Rem.
\\ \hline \hline
1&  & &${{\tau _0}^2}v$&$-{{1\over 2}{{\tau _0}^2}{v^2}}$&$
-{{1\over 6}{{\tau _0}^2}{v^3}}$& \\  \hline
2&  & &$-{{{\tau _0}^2}v}$&${1\over 2}{{\tau _0}^2}{v^2}$&${{1\over 6}
{\tau _0}^2}{v^3}$& \\ \hline
3&  &${{v_0}^2}{\tau _0}v$& & &${{v_0}^2}{\tau _0}a$& \\ \hline
4&  & & &${{\tau _0}^2}{v_0}v$&${1\over 2}{v_0}{{\tau _0}^2}{v^2}$&  \\
\hline
5&  & & & &$-{{v_0}^3}{\tau _0}\tau$ &  \\ \hline
6&  & & & &${{v_0}^2}{{\tau _0}^2}v$&  \\ \hline
7&  & & & &$-{{v_0}^2}{{\tau _0}^2}v$&  \\ \hline
8&  &${{v_0}^2}{\tau _0}v$& & &${{{v_0}^2}{\tau _0}a+{{v_0}^2}{{\tau _0}^2}v}
$& \\ \hline
9&  &${{v_0}^2}{\tau _0}v$& & &${{{v_0}^2}{\tau _0}a-{{v_0}^2}{{\tau _0}^2}v}
$& \\ \hline
10 & &${{v_0}^2}{\tau _0}v$&${{\tau _0}^2}v$&$-{{1\over 2}{{\tau _0}^2}
{v^2}}$&${{{v_0}^2}{\tau _0}a}-{{1\over 6}{{\tau _0}^2}{v^3}}$& \\  \hline
11 & &${{v_0}^2}{\tau _0}v$&$-{{\tau _0}^2}v$&${{1\over 2}{{\tau _0}^2}
{v^2}}$&${{{v_0}^2}{\tau _0}a}+{{1\over 6}{{\tau _0}^2}{v^3}}$& \\  \hline
12&  & & &${{\tau _0}^2}{v_0}v$&${1\over 2}{v_0}{{\tau _0}^2}{v^2}-
{{v_0}^3}{\tau _0}\tau$&  \\ \hline
13& ${\tau _0}{v_0}v$&$-{{1\over 2}{\tau _0}{v_0}{{v}^2}}$& & &$\varepsilon
{{v_0}^2}{{\tau _0}^2}v-{v_0}{\tau _0}m$&$\varepsilon \in$ R \\ \hline
14& ${\tau _0}{v_0}v$&$-{{1\over 2}{\tau _0}{v_0}{v^2}}$& &
${{\tau _0}^2}{v_0}v$&$\varepsilon{{v_0}^2}{{\tau _0}^2}v-{v_0}{\tau _0}m+
{1\over 2}{v_0}{{\tau _0}^2}{v^2}$&$\varepsilon \in$ R \\ \hline
15&  & &${{\tau _0}^2}v$&$-{{1\over 2}{{\tau _0}^2}
{v^2}}$&$-{{{v_0}^3}{\tau _0}\tau}-{1\over 6}{{\tau _0}^2}{v^3}$& \\  \hline
16&  & &$-{{\tau _0}^2}v$&${{1\over 2}{{\tau _0}^2}{v^2}}$
&$-{{{v_0}^3}{\tau _0}\tau}+{{1\over 6}{{\tau _0}^2}{v^3}}$& \\  \hline
17&  & & &$-{\tau _0}{{v_0}^2}\tau$&$-{{v_0}^2}{\tau _0}a$& \\ \hline
18&  & & &$-{\tau _0}{{v_0}^2}\tau$&$-{{\tau _0}{{v_0}^3}\tau}
-{{\tau _0}{{v_0}^2}a}$&   \\ \hline
19&  &${{v_0}^2}{\tau _0}v$& &$-\varepsilon{\tau _0}{{v_0}^2}\tau$
&$(1-\varepsilon){{v_0}^2}{\tau _0}a$&$\varepsilon \neq  0 $\\ \hline
20&  &${{v_0}^2}{\tau _0}v$& &${{\tau _0}{{v_0}^2}\tau}+{{{\tau _0}^2}{v_0}v}$
&$(1+1){{v_0}^2}{\tau _0}a+{1\over 2}{v_0}{{\tau _0}^2}{v^2}$&
 \\ \hline
21&  &$-{{v_0}^3}\tau$& &$-\varepsilon{\tau _0}{{v_0}^2}\tau$
&$-\varepsilon{{v_0}^2}{\tau _0}a-{1\over 2}{{\tau}^2}{{v_0}^3}$
&$\varepsilon\in $ R\\ \hline
22&  &$-{{v_0}^3}\tau$& &$-\varepsilon{\tau _0}{{v_0}^2}\tau$
&$-\varepsilon{{v_0}^2}{\tau _0}a-{1\over 2}{{\tau}^2}{{v_0}^3}+
{{v_0}^2}{{\tau _0}^2}v$&$\varepsilon\in $ R\\ \hline
23&  &$-{{v_0}^3}\tau$& &$-\varepsilon{\tau _0}{{v_0}^2}\tau$
&$-\varepsilon{{v_0}^2}{\tau _0}a-{1\over 2}{{\tau}^2}{{v_0}^3}-
{{v_0}^2}{{\tau _0}^2}v$&$\varepsilon\in $ R\\ \hline
24&  &$-{{v_0}^3}\tau$& &$-\varepsilon{\tau _0}{{v_0}^2}\tau+
{{\tau _0}^2}{v_0}v$
&$-\varepsilon{{v_0}^2}{\tau _0}a-{1\over 2}{{\tau}^2}{{v_0}^3}+{1\over 2}
{v_0}{{\tau _0}^2}{{v}^2}$&$\varepsilon\in $ R\\ \hline
25&  &$-{{v_0}^3}\tau$& &$-\varepsilon{\tau _0}{{v_0}^2}\tau-
{{\tau _0}^2}{v_0}{v}$
&$-\varepsilon{{v_0}^2}{\tau _0}a-{1\over 2}{{\tau}^2}{{v_0}^3}-{1\over 2}
{{v_0}}{{\tau _0}^2}{v^2}$&$\varepsilon\in $ R\\ \hline
26&  & &$-{v_0}{\tau _0}\tau$&${\tau _0}{v_0}a$&$-{v_0}{\tau _0}m$& \\ \hline
\end{tabular}
}
\end{center}

As far as Table 2 is concerned the following remark is in order. Up to
now we were dealing with dimensionless generators and group parameters.
In order to take care about proper dimensions we replace our generators
by dimensionful ones according to the rules:
\[ H\rightarrow {H\over {\tau _0}},\;P\rightarrow {P\over {{v_0}{\tau _0}}},\;
K\rightarrow {K\over {v_0}},\; M \rightarrow {M\over {{{v_0}^2}{\tau _0}}}\]
where $\tau _0$ and $v_0$ are arbitrary time and velocity units; the group
parameters are redefined appropriately. This redefinition has been already
taken into account in Table 2
\section{Conclusions}

We have classified all inequivalent bialgebra structures on the centrally
extended two-dimensional Galilei algebra and found the corresponding
Lie--Poisson structures on the group. The resulting classification appears
to be quite rich and contains 26 inequivalent cases, eigth of them being
the coboundary ones. This is in contrast with semisimple case as well as
the case of fourdimensional Poincare group where there are only coboundary
structures.
\section{Acknowledgment}
The Author acknowledges Prof.~P.~Kosi\'nski for a careful reading of the
manuscript and many helpful suggestion. Special thanks are also due to
Prof.~S.~Giller, Dr.~C.~Go\-ne\-ra, Prof.~P.~Ma\'slanka and MSc.~E.~Kowalczyk
for valuable discussion.


\end{document}